\documentclass[a4paper,11pt]{article}
\usepackage{graphicx}
\usepackage{amsfonts}
\usepackage{amsmath}
\usepackage[utf8]{inputenc}
\usepackage{url}
\usepackage{color}
 \usepackage{amssymb}
\usepackage{multirow}
\usepackage{overpic}

\newtheorem {Proposition}{Proposition}[section]
\newtheorem {Lemma}[Proposition] {Lemma}
\newtheorem {Theorem}[Proposition]{Theorem}
\newtheorem {Corollary}[Proposition]{Corollary}

\newtheorem {Example}{Example}[section]

\numberwithin{equation}{section}

\title{A Central Limit Theorem for $L_p$ transportation cost with applications to Fairness Assessment in  Machine Learning}

 \author{Eustasio del Barrio$^a$, Paula Gordaliza$^{a,b}$, and Jean-Michel Loubes$^b$\\
	$^{a}$\textit{IMUVA, Universidad de Valladolid}\\
	$^{b}$\textit{Institut de Math\'ematiques de Toulouse
}}

\begin{document}

\maketitle

\begin{abstract}
We provide a Central Limit Theorem for the Monge-Kantorovich distance between two empirical distributions with size $n$ and $m$, $W_p(P_n,Q_m)$ for $p>1$ for observations on the real line, using a minimal amount of assumptions. We provide an estimate of the asymptotic variance which enables to build a two sample test to assess the similarity between two distributions. This test is then used to provide a new criterion to assess the notion of fairness of a classification algorithm.
\end{abstract}

 \vskip .2in
\section{Introduction}

The analysis of the minimal transportation cost between two sets of random points or
of the transportation cost between an empirical and a reference measure is by now a classical
problem in probability, to which a significant amount of literature has been devoted. 
In the case of two sets of $n$ random points, say $X_1,\ldots,X_n$ and $Y_1,\ldots,Y_n$ in $\mathbb{R}^d$,
the object of interest is
$$T_{c,n}=\min_{\sigma}\frac 1 n\sum_{i=1}^n c(X_i,Y_{\sigma(i)}),$$
where $\sigma$ ranges is the set of permutations of $\{1,\ldots,n\}$ and $c(\cdot,\cdot)$ is some cost function. 
$T_{c,n}$ is usually referred to as the cost of optimal matching. This optimal matching problem is closely related
to the Kantorovich optimal transportation problem, which, in the Euclidean setting amounts to the minimization of
$$I[\pi]=\int_{\mathbb{R}^d\times \mathbb{R}^d} c(x,y) d\pi(x,y),$$
with $\pi$ ranging in the set of joint probabilities on $\mathbb{R}^d\times \mathbb{R}^d$
with marginals $P,Q$. Here $P$ and $Q$ are two probability measures on $\mathbb{R}^d$ and
the minimal value of $I[\pi]$ is known as the optimal transportation cost between $P$ and $Q$.
The cost function $c(x,y)=\|x-y\|^p$ has received special attention and we will write
$\mathcal{W}_p^p(P,Q)$ for the optimal transportation cost in that case. It is well known that
with this choice of cost function $T_{c,n}=\mathcal{W}_p^p(P_n,Q_n)$, with $P_n$ and $Q_n$ denoting
the empirical measures on $P_n$ and $Q_n$. 

How large is the cost of optimal matching, $\mathcal{W}_p^p(P_n,Q_n)$? Under the assumption that $X_1,\ldots,X_n$
are i.i.d. $P$, $Y_1,\ldots,Y_n$ are i.i.d. $Q$ and $P$ and $Q$ have finite $p$-th moment is is easy to conclude that
$\mathcal{W}_p^p(P_n,Q_n)\to \mathcal{W}_p^p(P,Q)$ almost surely. One might then wonder about the rate of approximation,
that is, how far is the empirical transportation cost from its theoretical counterpart.

\indent  Much effort has been devoted
to the case when $P=Q$, namely, when the two random samples come from the same random generator. In this case 
$\mathcal{W}_p^p(P,Q)=0$ and the goal is to determine how fast does the empirical optimal matching cost vanish.
It is known from the early works \cite{ajtai1984optimal} and \cite{Talagrand} that the answer depends on the dimension 
$d$. In the case when $P=Q$ is the uniform distribution on the unit hypercube $\mathcal{W}_p(P_n,Q_n)=O(n^{-1/d})$,
if $d\geq 3$, with a slightly worse rate if $d=2$. The results for $d\geq 3$ were later extended to a more general setup
in \cite{dobric1995asymptotics}, covering the case when $P=Q$ has bounded support and a density satisfying some smoothness
requirements. The one-dimensional case is different. If $p=1$ then, under some integrability assumptions
$\mathcal{W}_1(P_n,P)=O_P(n^{-1/2})$, with $\sqrt{n}\mathcal{W}_1(P_n,P)$ converging weakly to a non Gaussian limit, with, see \cite{del1999central}. 
If $p>1$ then it is still possible to get a limiting distribution for $\sqrt{n}\mathcal{W}_p(P_n,P)$, but now 
integrability assumptions are not enough and the available results require some smoothness conditions
on $P$ (and on its density), see \cite{del2005asymptotics} for the case $p=2$. In fact, see \cite{bobkov2014one}, the condition that $P$ has a positive density
an interval is necessary for boundedness of the sequence $\sqrt{n}E(\mathcal{W}_p(P_n,P))$ if $p>1$. \vskip .1in
Very recently a CLT in general dimension has been provided in \cite{del2017central}. The authors provide a CLT for quantities $\mathcal{W}_p^p(P_n,Q_n)$  concentrating around their mean $E (\mathcal{W}_p^p(P_n,Q_n))$ under some moment conditions (moments of of order $4+\delta$ with $\delta>0$ are required).  In this paper, we sharpen for the uni-dimensional case their results.   We prove asymptotic normality of $r_{n,m}(\mathcal{W}_p(P_n,Q_m)-\mathcal{W}_p(P,Q))$ for an increasing sequence $r_{n,m}$ under minimal moment and smoothness assumptions. \vskip .1in 

Such result enables to construct goodness of fit between two distributions  but also to assess  how similar two different
distributions $P$ and $Q$ can be. The similarity measure between the  distributions is the Wasserstein distance and we want to test if  $\mathcal{W}_p(P,Q)  \geq \Delta_0$ versus $\mathcal{W}_p(P,Q) < \Delta_0$ for $\Delta_0$ a chosen threshold. Note that in a different setting, this test is also considered in \cite{ramdas2017wasserstein}. \\
\indent  An application is given by  the recent  framework of Fair Learning or Disparate Impact Assessment which has received a growing attention driven by the generalization of machine learning in nowadays' life. We refer for instance to \cite{romei_ruggieri_2014},  \cite{pedreschi2012study}, \cite{2017arXiv170300056C} or \cite{2018arXiv180204422F} and references therein. In this setting, decisions are driven by machine learning procedures and the main concern is  to detect whether  a decision rule, learnt from variables $X$ is biased with respect to a subcategory of the population. For this a variable $S$ is denoted as a   protected attribute and splits the population into two groups $S=1$ and $S=0$. A decision rule is called unfair for $S$ when it exhibits a different behavior depending mainly on the values of $S$ and not on the values of the variables $X$. This discrimination may come from the algorithm or from a biased situation that would have been learnt from the learning sample. This framework has been proposed originally in \cite{FFMSV} and been further developed in \cite{del2018obtaining}.\\
\indent Many criteria have been given in the recent literature on Fair Learning to detect this situation (see in \cite{berk2017fairness} or \cite{myfairtest} for a review). A majority of these definitions consider that the decision should be independent from the protected attribute, which means that the decision should have similar behavior in both cases. Actually if $P_{|S=0}$ ( respectively $P_{|S=1}$) denotes the distribution of the classifier for $S=0$ (respectively $S=1$), then the complete fairness called Statistical Parity is obtained when these two distributions are the same, which corresponds to the independency of the decision with respect to the protected attribute. Therefore, the level of fairness could be quantified by estimating the similarity between $P_{|S=0}$ and $P_{|S=1}$. Hence we provide a new way of assessing fairness in machine learning by considering confidence intervals to test the similarity of these distributions with respect to Monge-Kantorovich distance. We study how this criterion behaves compared to standard criterion in the fair learning literature. \vskip .1in
The paper falls into the following parts. Section~\ref{s:main} provides the main result, i.e the Central Limit Theorem for $L_p$ transportation cost for $p>1$. Section~\ref{s:simul} is devoted to some simulations while Section~\ref{s:fair} is devoted to the application of this test to detect disparate impact. Proofs are gathered in the Appendix.

\section{CLT for $L_p$ transportation cost} \label{s:main}

In this section we present the main results in this paper, namely, CLT's for the 
transportation cost between an empirical measure and a target measure or between
two empirical measures. 

\bigskip

To present our results, we set $h_p(x)=|x|^p$, $x\in \mathbb{R}$, $p>1$
and consider the functions
\begin{equation}\label{GausspCost}
c_p(t;F,G):=\int_{F^{-1}(\frac 1 2)}^{F^{-1}(t)} h_p'\big(s-G^{-1}(F(s))\big) ds,\quad 0<t<1.
\end{equation}
We note that $h'_p(x)=p\,\mbox{sgn}(x)|x|^{p-1}$. Since $F^{-1}(\frac 1 2)\leq s<
F^{-1}(t)$ implies $\frac 1 2\leq F(s)<t$ while for $F^{-1}(t)\leq s<
F^{-1}(\frac 1 2)$ we have $t\leq F(s) <\frac 1 2$, we see that $c_p(t;F,G)$ is finite for every
$t\in (0,1)$. Under the assumption $F,G\in\mathcal{F}_{2p}$ we show in Lemma \ref{integrability} below 
that, in fact, $c_p(\cdot;F,G)\in L_2(0,1)$. This allows us to introduce also
\begin{equation}\label{GausspCost2}
\bar{c}_p(t;F,G):=c_p(t;F,G)-\int_0^1 {c}_p(s;F,G)ds ,\quad 0<t<1.
\end{equation}
We observe that changing $F^{-1}(\frac 1 2)$ by $F^{-1}(t_0)$ in (\ref{GausspCost})
would not affect the definition of $\bar{c}_p(\cdot;F,G)$.

\begin{Lemma}\label{integrability}
	If $F, G\in \mathcal{F}_{2p}(\mathbb{R})$, $p>1$, then $c_p(\cdot;F,G)\in L_2(0,1)$ and $\bar{c}_p(\cdot;F,G)\in L_2(0,1)$.
	Furthermore, if $F_m, G_m\in\mathcal{F}_{2p}(\mathbb{R})$
	satisfy $\mathcal{W}_{2p}(F_m,F)\to 0$, $\mathcal{W}_{2p}(G_m,G)\to 0$ and $G^{-1}$ is continuous on $(0,1)$ then 
	$\bar{c}_p(\cdot;F_m,G_m)\to 
	\bar{c}_p(\cdot;F,G)$ in $L_2(0,1)$ as $m\to\infty$.
\end{Lemma}

\medskip
\noindent \textbf{Proof.} We set
$d_p=\max(1,2^{p-2})$, $p>1$, and observe that
\begin{eqnarray}\label{cotaintunif}
|c_p(t;F,G)|&\leq & p d_p \Big|\int_{F^{-1}(\frac 1 2)}^{F^{-1}(t)} \big(|s|^{p-1}+ |G^{-1}(F(s))|^{p-1}\big)dt\Big|\\
\nonumber
&\leq & d_p |F^{-1}(t)-F^{-1}({\textstyle \frac 1 2})|\Big(|F^{-1}(t)|^{p-1}+|F^{-1}({\textstyle \frac 1 2})|^{p-1}
+|G^{-1}(t)|^{p-1}+|G^{-1}({\textstyle \frac 1 2})|^{p-1} \Big).
\end{eqnarray}
The first claim follows upon using H\"older's inequality to check that $\int_0^{1} |F^{-1}(s)|^{2}|G^{-1}(s)|^{2(p-1)} ds<\infty$.
For the second observe that $\mathcal{W}_p(F_m,F)\to 0$ implies that $F_m^{-1}(t)\to F^{-1}(t)$ for every $t$ of continuity
for $F^{-1}$ (hence, for almost every $t\in (0,1)$) and also that $|F_m^{-1}|^{2p}$ is uniformly integrable (and the same holds
for $G_m^{-1}$, with convergence at every point in $(0,1)$ since $G^{-1}$ is continuous). By the remarks before
this Lemma we can assume without loss of generality that $F^{-1}$ is continuous at $\frac 1 2$. Then 
$|c_p(t;F_m,G)|\to |c_p(t;F,G)|$ at every $t$ of continuity for $F^{-1}$. Using the bound (\ref{cotaintunif}) for $c_p(t;F_m,G)$ and
the fact that $|F_m^{-1}|^{2p}$ and $|G_m^{-1}|^{2p}$ are uniformly integrable, we see that the sequence
$c_p^2(\cdot;F_m,G_m)$ is uniformly integrable and conclude that ${c}_p(\cdot;F_m,G)\to {c}_p(\cdot;F,G)$ and
$\bar{c}_p(\cdot;F_m,G_m)\to \bar{c}_p(\cdot;F,G)$ in $L_2(0,1)$. 
\quad $\Box$

\medskip
It is convenient at this point to introduce the notation 
\begin{equation}\label{VarF}
\sigma^2_p(F,G)=\int_0^1 \bar{c}^2_p(t;F,G) dt.
\end{equation}
Lemma \ref{integrability} ensures that $\sigma^2_p(F,G)$ is a finite constant provided $F$ and $G$ have finite moments of order $2p$. 
Also, if $F\ne G$ then $G^{-1}\circ F$, which is the optimal transportation map from $F$ to $G$ is different from the identity on a set of positive measure
and $\sigma^2_p(F,G)>0$ if $F$ is not a Dirac measure. We remark that $\sigma^2_p(F,G)$ is not, in general, symmetric in $F$ and $G$. 

\medskip
We are ready now for the main result in this section.
\begin{Theorem}\label{CLTLpCost}
	Assume that $F,G\in\mathcal{F}_{2p}(\mathbb{R})$ and $G^{-1}$ is continuous on $(0,1)$. 
	Then 
	\begin{itemize}
		\item[(i)] If $X_1,\ldots,X_n$ are i.i.d. $F$ and $F_n$ is the empirical d.f. based on the $X_i$'s
		$$\sqrt{n}(\mathcal{W}_p^p(F_n,G)-E\mathcal{W}_p^p(F_n,G))\to_w N(0,\sigma^2(F,G)).$$
		\item[(ii)] If, furthermore, $F^{-1}$ is continuous, $Y_1,\ldots,Y_m$ are i.i.d. $G$, independent of the $X_i$'s,
		$G_m$ is the empirical d.f. based on the $Y_j$'s and $\frac{n}{n+m}\to\lambda \in (0,1)$ then
		$$\textstyle \sqrt{\frac{nm}{n+m}}(\mathcal{W}_p^p(F_n,G_m)-E\mathcal{W}_p^p(F_n,G_m))\to_w N(0,(1-\lambda)\sigma^2(F,G)+\lambda\sigma^2(G,F) ).$$
	\end{itemize}
\end{Theorem}

A proof of this result is given in the Appendix. We would like to make some remarks about Theorem \ref{CLTLpCost} at this point.
There has been a significant interest in empirical transportation costs in recent times in the literature.
We should mention at least \cite{fournier2015rate}, giving moment bounds and concentration results for empirical
transportation with $L_p$ cost in general dimension, and \cite{bobkov2014one}, with a comprehensive discussion
of the one dimensional case. Both papers focus on the case where the law underlying the empirical measure 
and the target measure are equal (in the setup of Theorem \ref{CLTLpCost}, the case $F=G$). With the more specific
goal of CLT's for empirical transportation costs, \cite{sommerfeld2018inference} considers the case when the underlying 
probabilities are finitely supported, while \cite{tameling2017empirical} covers probabilities with countable 
support. The approach in these two cases relies on Hadamard directional differentiability of the dual form
of the finite (or countable) linear program associated to optimal transportation. Without the constraint of
countable support, \cite{del2017central} covers quadratic transportation costs in general dimension.

There are
similarities between the approach in \cite{del2017central} and the presentation here, as one can see from a look
at our Appendix. We must emphasize some significant differences, however. 
An obvious one is that here we only deal with one dimensional probabilities. On the other hand, we 
cover general $L_p$ costs. A more significant difference is that assumptions in Theorem \ref{CLTLpCost} 
are sharp. Let us focus on \textit{(i)} to discuss this point. To make sense of $\mathcal{W}_p^p(F_n,G)$
we must consider $G$ with finite $p$-th moment. Now, if we want $F$ to satisfy \textit{(i)} 
for every $G$ with finite $p$-th moment, by taking $G$ to be Dirac's measure on $0$ we see that 
$F$ must have finite moment of order $2p$. Then it is easy to check that, $\sigma^2(F,G)<\infty$ for all $F$ 
with finite moment of order $2p$ if and only if $G$ has a finite moment of order $2p$. Thus, the 
assumption of finite moments of order $2p$ for $F$ and $G$ seems to be a minimal requirement for 
\textit{(i)} to hold. We note that for the quadratic cost, $p=2$, Theorem 4.1 in \cite{del2017central} 
required finite moments of
order $4+\delta$ on $P$ and $Q$ for some $\delta>0$.

Some words on the role of the continuity of $G^{-1}$ in \textit{(i)} are also in place here. 
That some sort of regularity of the quantile function is needed for handling the empirical transportation
functional in dimension one was observed in \cite{bobkov2014one}. In the case $F=G$, absolute continuity
of $F^{-1}$ is a necessary condition for having $E(\mathcal{W}_p(F_n,F))=O(\frac{1}{\sqrt{n}})$ 
(Theorem 5.6 in \cite{bobkov2014one}). Continuity of $G^{-1}$ is also related to assumption (3) in 
\cite{del2017central}. In fact, that assumption, in the case of one-dimensional probabilities, implies
that $G$ is supported in a (possibly unbounded) interval and $G^{-1}$ is differentiable in 
the interior of that interval. Hence, the regularity assumption in Theorem \ref{CLTLpCost} is also 
slightly weaker that that in Theorem 4.1 in \cite{del2017central}. We should also note at this point that
Theorem 1 in \cite{sommerfeld2018inference}, for the case finitely supported probabiities on the real line  
corresponds to a case of discontinuity of the quantile functions and this can lead to nonnormal limiting
distributions.

\medskip
We would also like to discuss the role of the centering constants in Theorem \ref{CLTLpCost}. 
Under more restrictive assumptions there are similar CLT's in which $E\mathcal{W}_p^p(F_n,G)$ is replaced
by the simpler constants $\mathcal{W}_p(F,G)$ (see, e.g., Therem 4.3 in \cite{del2017central}).
In fact, the Kantorovich duality (see, e.g., \cite{Villani2003})
yields that
$$\mathcal{W}_p^p(F,G)=\sup_{(\varphi,\psi)\in\Phi_p} \int \varphi dF+\int \psi dG,$$
where $\Phi_p$ is the set of pairs of integrable functions (with respec to $F$ and $G$, respectively)
satisfying $\varphi(x)+\psi(y)\leq |x-y|^p$. But this entails $E(\mathcal{W}_p^p(F_n,G))\geq
\sup_{(\varphi,\psi)\in\Phi_p} E\big(\int \varphi dF_n\big)+\int \psi dG=
\sup_{(\varphi,\psi)\in\Phi_p} \int \varphi dF+\int \psi dG=\mathcal{W}_p^p(F,G)$. Hence, we can replace 
the centering constants in Theorem \ref{CLTLpCost} provided
\begin{equation}\label{centeringprop}
0\leq \sqrt{n}\big(E(\mathcal{W}_p^p(F_n,G))-\mathcal{W}_p^p(F,G)\big)\to 0.
\end{equation}

Finding sharp conditions under which (\ref{centeringprop}) holds seems to be a delicate issue.
We limit ourselves to providing a set of sufficient conditions for it.
The case $F=G$ has been considered in \cite{bobkov2014one} and can be handled with simple moment
conditions. The general case that we consider here seems to add some smoothness requirements.
We limit our discussion to $p\geq 2$.
We will assume that $F$ is twice differentiable, with nonvanishing density, $f$, in the interior 
of $\mbox{supp}(F)=\mbox{cl}\{x:\, F(x)\notin \{0,1\}\}$ and satisfies
\begin{equation}\label{smoothnessplus}
\sup_{t\in (0,1)} \frac{t(1-t)|f'(F^{-1}(t))|}{f^2(F^{-1}(t))}<\infty.
\end{equation}
Furthermore, we will assume that
\begin{equation}\label{tech1}
\mbox{for some } s\in \textstyle (\frac p 4, \frac p 2), \quad n^s E\mathcal{W}_p^p(F_n,F)\to 0\quad \mbox{ as } n\to\infty, 
\end{equation}
\begin{equation}\label{tech2}
\frac{1}{\sqrt{n}}\int_{\frac 1 n}^{1-\frac 1 n} \frac{(t(1-t))^{1/2}}{f^2(F^{-1}(t))}dt\to 0,
\end{equation}
\begin{equation}\label{tech3}
\int_0^1\int_0^1 \frac{(s\wedge t-st)^2}{f^2(F^{-1}(s))f^2(F^{-1}(t))}dsdt<\infty.
\end{equation}
Condition (\ref{smoothnessplus}) is a natural condition for approximating the quantile process by
a weighted uniform standard process. We refer to \cite{del2005asymptotics} for details. The other three conditions are implied
by the stronger assumption
\begin{equation}\label{tech4}
\int_{0}^{1} \frac{(t(1-t))^{p/2}}{f^p(F^{-1}(t))}dt<\infty.
\end{equation}
This condition is, essentially, needed for ensuring that $n^{p/2} E\mathcal{W}_p^p(F_n,F)$ is a bounded sequence,
see \cite{bobkov2014one}. We would like to note that (\ref{tech4}) does not hold for Gaussian $F$, while
(\ref{tech1}), (\ref{tech2}) and (\ref{tech3}) do.

With these assumptions we can prove the following.
\begin{Proposition}\label{Proposition}
	Assume $p\geq 2$. Under the assumptions of Theorem \ref{CLTLpCost},
	\begin{itemize}
		\item[(i)] if $F$ satisfies (\ref{smoothnessplus}) to (\ref{tech3}) then (\ref{centeringprop}) holds and, 
		as a consequence,
		$$\sqrt{n}(\mathcal{W}_p^p(F_n,G)-\mathcal{W}_p^p(F,G))\to_w N(0,\sigma^2(F,G)).$$
		\item[(ii)] if, furthermore, $G$ satisfies (\ref{smoothnessplus}) to (\ref{tech3}) then
		$$\textstyle \sqrt{\frac{nm}{n+m}}(\mathcal{W}_p^p(F_n,G_m)-E\mathcal{W}_p^p(F,G))\to_w N(0,(1-\lambda)\sigma^2(F,G)+\lambda\sigma^2(G,F) ).$$
	\end{itemize}
\end{Proposition}

A proof of Proposition \ref{Proposition} is given in the Appendix. The scheme of proof, in fact, relies
on some auxiliary results in \cite{del2005asymptotics} that give, through a completely different approach, 
asymptotic normality of $\sqrt{n}(\mathcal{W}_p^p(F_n,G)-\mathcal{W}_p^p(F,G))$. \vskip .1in
The economy in assumptions that 
one can gain from dealing with the centering in Theorem \ref{CLTLpCost} is, in our view, remarkable.
Providing sharper conditions under which (\ref{centeringprop}) holds remains an interesting open question.

\bigskip
For the statistical application of Theorem \ref{CLTLpCost} it is of interest to have a 
consistent estimator of the asymptotic variances. In the two sample case 
this can be done as follows. Define
$$d_{i,n,m}(X,Y)=\sum_{j=2}^i \Big[ \big|X_{(j)}-G_{m}^{-1}(\textstyle \frac {j-1}n)\big|^p-
\big|X_{(j-1)}-G_{m}^{-1}(\textstyle \frac {j-1}n)\big|^p\Big],\quad i=2,\ldots,n$$
with $d_{1,n,m}(X,Y)=0$ and
\begin{equation}\label{varest}
\hat{\sigma}_{1,n,m}^2=\textstyle \frac 1 n\sum_{i=1}^n d^2_{i,n,m}(X,Y)-\Big(\textstyle \frac 1 n\sum_{i=1}^n d_{i,n,m}(X,Y)\Big)^2.
\end{equation}
We define $\hat{\sigma}_{2,n,m}^2$ similarly exchanging the roles of the $X_i$'s and the $Y_j$'s. Finally, we set
\begin{equation}\label{varest2}
\hat{\sigma}_{n,m}^2=\textstyle \frac m{n+m}\hat{\sigma}_{1,n,m}^2+\frac n{n+m}\hat{\sigma}_{2,n,m}^2.
\end{equation}
We show next that $\hat{\sigma}_{n,m}^2$ is a consistent estimator of the asymptotic variance in the two sample case in Theorem \ref{CLTLpCost}.
A consistent estimator for the asymptotic variance in the one sample case can be obtained similarly. We omit details.

\begin{Proposition}\label{consistentsigma}
	If $F,G\in\mathcal{F}_{2p}(\mathbb{R})$ and $F^{-1},G^{-1}$ are continuous on $(0,1)$ then 
	$$\hat{\sigma}_{n,m}^2\to (1-\lambda)\sigma^2(F,G)+\lambda\sigma^2(G,F)$$
	almost surely.
\end{Proposition}
\medskip
\noindent \textbf{Proof.} Simply note that $\hat{\sigma}_{1,n,m}^2=\int_0^1 \bar{c}_p^2(t;F_n,G_m)dt$ and apply Lemma
\ref{integrability}. \hfill $\Box$

\bigskip
As a consequence of Propositions \ref{Proposition} and \ref{consistentsigma} we have that if, additionally,
$$F\ne G$$ and $F$ (or $G$) is not a Dirac measure then
\begin{equation}\label{CLTstandard}
\textstyle \sqrt{\frac{nm}{n+m}}\frac{(\mathcal{W}_p^p(F_n,G_m)-\mathcal{W}_p^p(F,G))}{\hat{\sigma}_{n,m}}\to_w N(0,1).
\end{equation}
We can use (\ref{CLTstandard}) for statistical applications in several ways. From (\ref{CLTstandard}) we see that 
\begin{equation}\label{asymptCI}
\big[\mathcal{W}_p^p(F_n,G_m)\pm \textstyle \sqrt{\frac{n+m}{nm}} \hat{\sigma}_{n,m} \Phi^{-1}(1-\frac \alpha 2)\big]
\end{equation}
is a confidence interval for $\mathcal{W}_p^p(F,G)$ with asymptotic confidence level $1-\alpha$. Alternatively, we could
consider the testing problem
\begin{equation}\label{testingproblem}
H_0:\, \mathcal{W}_p(F,G)\geq \Delta_0,\quad \mbox{ vs } \quad H_1:\,  \mathcal{W}_p(F,G)< \Delta_0,
\end{equation}
where $\Delta_0$ is some threshold (to be determined by the practitioner). Rejection of the null in (\ref{testingproblem})
would yield statistical evidence that the d.f.'s $F$ and $G$ are almost equal. We can handle this problem by rejecting the null if
\begin{equation}\label{rejectionregion}
\mathcal{W}_p^p(F_n,G_m)<\Delta_0^p-\textstyle\sqrt{\frac{n+m}{nm}} \hat{\sigma}_{n,m} \Phi^{-1}(1- \alpha).
\end{equation}
It follows from (\ref{CLTstandard}) that the test defined by (\ref{rejectionregion}) has asymptotic level $\alpha$. In the next 
section we explore the use of this test for the assessment of fairness of learning algorithms.

\section{Simulations and Results} \label{s:simul}
In this section, we first analize the consistency of the variance estimation given by \eqref{varest}-\eqref{varest2} established in Proposition \ref{consistentsigma}. Then, we check the performance of the test and finally, we apply both tools to the Fair Learning problem.

Consider two independent samples $X_1,\ldots,X_n$ i.i.d. and $Y_1,\ldots,Y_m$ i.i.d. of distributions $F$ and $G$, respectively, and denote by $F_n$ and $G_m$ the corresponding empirical distribution function on each sample. We have simulated these samples undergoing the following models, for which we can compute the exact expression for the asymptotic variance in Proposition \ref{consistentsigma}.
\begin{Example}[Location model]\label{example:location}
	 Consider $F\sim N(0,1)$ and $G \sim N(\mu,1), \ \mu \in \mathbb{R}$. We can write $G^{-1}(t)=\Phi^{-1}(t)+\mu$ and the Wasserstein distance between both distributions is $W_p(F,G)=\left|\mu\right|, \ p \geq 1$. In this case, we can compute the functions
	\begin{align*}
	c_p(t;F,G)&=-p.\text{sgn}(\mu)\left|\mu\right|^{p-1}\Phi^{-1}(t)\\
	c_p(t;G,F)&=p.\text{sgn}(\mu)\left|\mu\right|^{p-1}\Phi^{-1}(t),
	\end{align*}
	and then
	\begin{align*}
	&\int_{0}^1c_p(t;F,G)^2dt=p^2\mu^{2(p-1)}\int_{0}^1(\phi^{-1}(t))^2dt=p^2\mu^{2(p-1)}\\
	&\int_{0}^1c_p(t;F,G)dt=p\left|\mu\right|^{p-1}\int_{0}^1\phi^{-1}(t)dt=0.
	\end{align*}
	Hence, in this model we have an expression for the true variance $\sigma^2=(1-\lambda)\sigma^2(F,G) +\lambda\sigma^2(G,F),$ where $\sigma(F,G)=\sigma(G,F)=p\mu^{p-1}$. In Table \ref{tab:varianceestimation}, we can see the estimation for an increasing size $n=m$ of the samples, which are close in the limit to the true value.
\end{Example}

\begin{Example}[Scale-location model]\label{example:scale-location}
	Consider $F\sim N(0,1)$ and $G \sim N(\mu,\lambda), \ \left(\mu, \lambda\right) \in \mathbb{R}\times\mathbb{R}^+$. In this case, $G^{-1}(t)=\lambda\Phi^{-1}(t)+\mu$, and then
	\begin{align*}
	c_p(t;F,G)&=\frac{1}{1-\lambda}\left[\left|(1-\lambda)\phi^{-1}(t)-\mu\right|^p-\left|\mu\right|^p\right]\\
	c_p(t;G,F)&=\frac{\lambda}{\lambda-1}\left[\left|(\lambda-1)\phi^{-1}(t)+\mu\right|^p-\left|\mu\right|^p\right].
	\end{align*}
	
\end{Example}

\begin{table}[h]
	\centering
	\begin{tabular}{|c|c|c|c|}
		\hline
		n&p=1&p=2&p=3\\
		\hline
		50&0.7811&5.9767&9.4721\\
		
		100&0.8742&3.0618&9.668\\
		
		200&0.9262&5.1305&10.1345\\
		
		400&1.0510&4.9746&8.7785\\
		
		500&1.0023&4.0164&9.2851\\
		
		800&0.9858&3.4522&8.592\\
		
		1000&1.0923&4.399&8.9125\\
		
		2000&0.9868&3.4341&9.1057\\
				
		5000&0.9932&4.1488&8.9690\\
		
		10000&0.9999&4.0661&9.1961\\
				
		20000&0.9842&4.0426&8.9744\\
		
		50000&1.003&3.9567&9.1324\\
		
		100000&0.9965&4.0184&8.9922\\
		\hline
		$\sigma^2$&1&4&9\\
		\hline
	\end{tabular}
	\caption{Estimates of the variance of the asympotic distribution in the location model of Example \ref{example:location} with $\mu=1$}
	\label{tab:varianceestimation}
\end{table}

\begin{table}[h]
	\centering
	\begin{tabular}{|c|c|c|c|c|c|}
		\hline
		$p$&n&$\mu$=1&$\mu$=0.9&$\mu$=0.7&$\mu$=0.5\\
		\hline
		\multirow{8}{*}{1} &
		50&0.062&0.146&0.481&0.825\\
		
		&100&0.055&0.193&0.698&0.974\\
		
		&200&0.053&0.275&0.918&1\\
		
		&400&0.051&0.413&0.995&1\\
		
		&500&0.051&0.481&0.999&1\\
		
		&800&0.052&0.64&1&1\\
		
		&1000&0.054&0.728&1&1\\
		
		&2000&0.047&0.937&1&1\\
		\hline
		\multirow{8}{*}{2} &
		50&0.074&0.167&0.513&0.839\\
		
		&100&0.063&0.198&0.717&0.979\\
		
		&200&0.059&0.272&0.927&1\\
		
		&400&0.055&0.422&0.995&1\\
		
		&500&0.05&0.484&0.999&1\\
		
		&800&0.053&0.651&1&1\\
		
		&1000&0.053&0.736&1&1\\
		
		&2000&0.051&0.935&1&1\\
		\hline
		\multirow{8}{*}{3} &
		50&0.071&0.154&0.515&0.822\\
		
		&100&0.0662&0.206&0.715&0.973\\
		
		&200&0.057&0.266&0.925&1\\
		
		&400&0.052&0.422&0.992&1\\
		
		&500&0.057&0.497&0.997&1\\
		
		&800&0.053&0.652&1&1\\
		
		&1000&0.053&0.733&1&1\\
		
		&2000&0.051&0.937&1&1\\
		\hline
	\end{tabular}
	\caption{Estimated probabilities of rejection in the location model with $\Delta_0=1$}
	\label{tab:alphalocationmodel}
\end{table}

%
%
%
%
%
%
%
%
%
%
%
%
%
%
%
%
%
%
%
%
%

\begin{table}[h]
	\centering
	\begin{tabular}{|c|c|c|c|c|c|c|}
		\hline
		$\Delta_0$&$p$&n&$\mu$=1, $\lambda$=2&$\mu$=1, $\lambda$=1.5&$\mu$=0, $\lambda$=2&$\mu$=0, $\lambda$=1.5\\
		\hline
		\multirow{8}{*}{$W_1(N(0,1),N(1,2))=0.7978846$} & \multirow{8}{*}{1} &
		50&0.047&0.165&0.535&0.996\\
		
		&&100&0.045&0.195&0.8&1\\
		
		&&200&0.036&0.323&0.974&1\\
		
		&&400&0.052&0.532&1&1\\
		
		&&500&0.056&0.614&1&1\\
		
		&&800&0.035&0.810&1&1\\
		
		&&1000&0.045&0.895&1&1\\
		
		&&2000&0.050&0.994&1&1\\
		\hline
		\multirow{8}{*}{$W_2(N(0,1),N(1,2))=2$} & \multirow{8}{*}{2} &
		50&0.078&0.376&0.595&0.998\\
		
		&&100&0.067&0.551&0.823&1\\
		
		&&200&0.062&0.786&0.976&1\\
		
		&&400&0.055&0.969&1&1\\
		
		&&500&0.059&0.985&1&1\\
		
		&&800&0.052&1&1&1\\
		
		&&1000&0.056&1&1&1\\
		
		&&2000&0.05&1&1&1\\
		\hline
		\multirow{8}{*}{$W_3(N(0,1),N(1,2))=1.611195$} & \multirow{8}{*}{3} &
		50&0.091&0.569&0.571&0.997\\
		
		&&100&0.093&0.762&0.758&1\\
		
		&&200&0.072&0.935&0.939&1\\
		
		&&400&0.06&1&0.996&1\\
		
		&&500&0.064&0.999&0.997&1\\
		
		&&800&0.069&1&1&1\\
		
		&&1000&0.06&1&1&1\\
		
		&&2000&0.049&1&1&1\\
		\hline
	\end{tabular}
	\caption{Estimated probabilities of rejection in the scale-location model}
	\label{tab:alphascalelocationmodel}
\end{table}
To check the performance of the test \eqref{testingproblem}, we have simulated observations in the scenarios of the Examples \ref{example:location} and \ref{example:scale-location} for different values of the parameters of location and scale. Table \ref{tab:alphalocationmodel} shows the estimated frequencies of rejection of the test in the location model when $\Delta_0=1$, and different values for the cost $p=1,2,3$.\vskip .1in Under the null hypothesis $H_0$, that is when $\mu=1$, we see that the covering level achieves the nominal value $\alpha=0.05$. Moreover, under the alternative $H_1$, that is when $\mu=0.5, 0.7, 0.9$, the values show that the test has high power. Similar results are obtained for the scale-location model, which are contained in Table \ref{tab:alphascalelocationmodel}, for different values of the threshold $\Delta_0$ and the cost $p=1,2,3$. \\ \indent Under the null hypothesis $H_0$, that is when $\left(\mu, \lambda\right)=(1,2)$, the estimated level reaches the nominal value $\alpha=0.05$. The test shows again high power when $\left(\mu, \lambda\right)=(1,1.5), (0,2), (0,1.5)$. We note that even in the case $p=1$, which is not covered by the theoretical results in this paper, the simulations in both models show that the test has asymptotically level $\alpha$ and that its power is very high is most cases.

\section{Application to Fair Learning} \label{s:fair}
Fair learning is devoted to the analysis of bias that appear when learning automatic decisions (mainly classification rules) from a learning sample . This sample may be prejudiced against a population, which means that the variable to be predicted is, in the sample, unbalanced between the two groups. Hence when trying to find a classification rule, the algorithm will use the discrimination present in the sample rather than learning a true link function.  This bias can have been set intentionally or may reflect the bias present in the use cases. A striking example is provided by looking at banks predicting high income from a set of parameters in order to grant a loan. Despite of their claim, this prediction leads to a clear discrimination between male and female, while this variable should not play any role in such forecast. Yet the constitution of the learning sample lead the classifier to underestimate the income of female compared to male. \\
\indent Hence it is important to detect such automatic bias in order to prevent a generalization and even worse a justification of a discriminatory behavior. Many criterion have been proposed to quantify the influence of the group variable $S$ in the behavior of the machine learning algorithm, most of them consider a notion of similarity between the decisions of algorithm conditionally to the belonging of each group. In the following we propose a new criterion by considering the Monge-Kantorovich distance between the distribution of the classification rule conditioned by the two groups. This problem is at the heart of recent studies in machine learning, leading to a new field of research called fair machine learning. Hence To illustrate the application of the theory in previous section to the problem of fairness in Machine Learning, we consider the \textit{Adult Income} data set (available at \url{https://archive.ics.uci.edu/ml/datasets/adult}). It contains $29.825$ instances consisting in the values of $14$ attributes, $6$ numeric and $8$ categorical, and a categorization of each person as having an income of more or less than $50,000\$$ per year. \vskip .1in

Recently, in \cite{del2018obtaining} the problem of forecasting a binary variable $Y \in \{0,1\}$ using observed covariates $X \in \mathbb{R}^d, \ d\geq 1,$ and assuming that the population is divided into two categories that represent a bias, modeled by a protected variable $S \in \{0,1\}$, is considered. The criteria of fairness in classification problems considered are called Disparate Impact (DI) and Balanced Error Rate (BER), which were introduced in \cite{FFMSV}.  For a classification rule $g$, Disparate Impact is a score that measures how the two probability $P(g(X)=1|S=0)$ and $P(g(X)=1|S=1)$ are close. The Balance Error Rate describes how the variable $S$ can be learnt by the classification rule $g$ originally meant to predict the variable $Y$. Using these criteria, they designed procedures to remove the possible discrimination, both in a partial and total way, that are based on the idea of moving the distributions of the variable $X$ conditionally given the value of the protected $S$. This approach originally proposed in  \cite{FFMSV} can also be found in \cite{DBLP:journals/corr/abs-1712-07924}.

Also in \cite{del2018obtaining}, confidence intervals for the empirical counterpart of DI proposed in \cite{myfairtest} and some numerical results of the procedures that remove discrimination are given for the Adult Income Data and the logit classifier. This classifier is used to make the prediction using the five numerical variables: \textit{Age, Education Level, Capital Gain, Capital Loss} and \textit{Worked hours per week}. Among the rest of the categorical attributes, the sensitive attribute to be the potentially protected is the \textit{Gender} (\textquotedblleft male\textquotedblright or \textquotedblleft female\textquotedblright). While in that work the logit is used for binary classification whether a person earns more or less than $50,000\$$ per year, here we will consider the result of the logistic regression, that is, the estimated  probability of positive outcome, which provides for all observation a distribution on the real line.  \\
This estimation is used to predict whether an individual will have a high income. This forecast algorithm presents some bias with respect to the gender in the sense that the learning sample is biased such that female with similar characteristics as male are more unlikely to be predicted that they will get a high income. This unfairness is usually shown using disparate impact assessment as discussed in \cite{del2018obtaining}. \vskip .1in

\begin{figure}[h]
	\centering
	\begin{overpic}[scale=.4,unit=1mm]{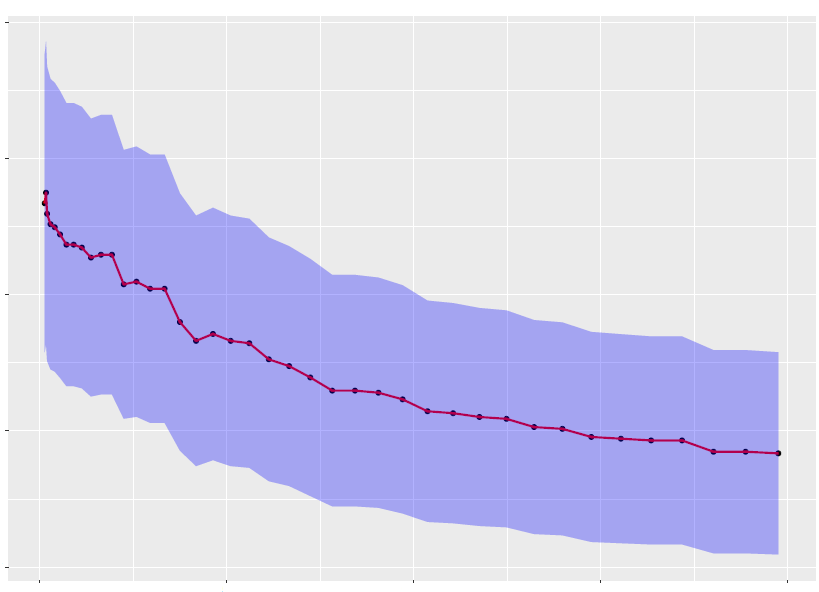}
		\put (-15,30){$\hat{DI}$}
		\put (-6,68){$1.2$}
		\put (-5,52){$1$}
		\put (-6,35){$0.8$}
		\put (-6,18){$0.6$}
		\put (-6,2){$0.4$}
		\put (50,-10){$W_2^2$}
		\put (4,-3){$0$}
		\put (23,-3){$0.003$}
		\put (46,-3){$0.006$}
		\put (70,-3){$0.009$}
		\put (92,-3){$0.012$}
		
	\end{overpic}
\vspace{1cm}
	\caption{Relationship between DI and $W_2^2(\mathcal{L}(g(X)\mid S=0),\mathcal{L}(g(X)\mid S=1))$}
	\label{fig:DI_W2}
\end{figure}

\begin{figure}[h]
	\centering
	\begin{overpic}[scale=.4,unit=1mm]{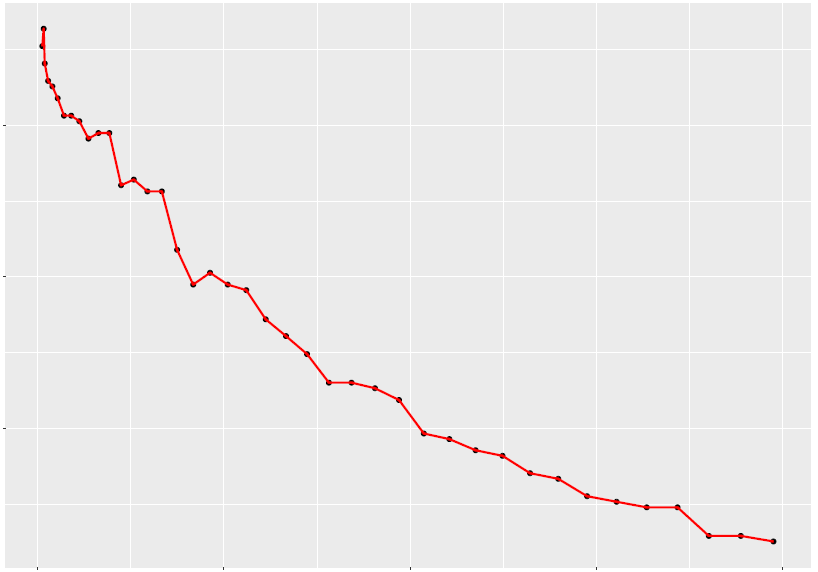}
	\put (-15,30){$\hat{BER}$}

	\put (-9,53){$0.49$}
	\put (-9,35){$0.48$}
	\put (-9,16){$0.47$}

	\put (50,-10){$W_2^2$}
	\put (4,-4){$0$}
	\put (23,-4){$0.003$}
	\put (46,-4){$0.006$}
	\put (70,-4){$0.009$}
	\put (92,-4){$0.012$}
	
\end{overpic}
\vspace{1cm}
\caption{Relationship between BER and $W_2^2(\mathcal{L}(g(X)\mid S=0),\mathcal{L}(g(X)\mid S=1))$}
\label{fig:BER_W2}
\end{figure}

We want to see if the test \eqref{testingproblem} is an appropriate tool to asses fairness in algorithmic classification results and could replace the Disparate impact.\\
\indent  Actually, fairness should be achieved as soon as the distribution of the forecast probability is close for the two groups corresponding to $S=0$ and $S=1$. We choose to directly control this closeness using Wasserstein distance and not using the Disparate Impact criterion. Yet, we study the relationship between the notions of Disparate Impact and the Balanced Error Rate of a classifier $g$, with the Wasserstein distance between the probability distributions of $g(X)=\log(1/(1+\exp(-\beta^T X))$ conditionally given the protected attribute $S=0,1$. In \cite{del2018obtaining}, it is proved that the BER is related to the distance in Total variation between these two conditional distributions, but that the optimal transportation cost is still a reasonable way of quantifying the distances. \\
\indent Hence we  study on simulations how the variation in the Wasserstein distance between the two distributions $\mathcal{L}(g(X)\mid S=0)$ and $\mathcal{L}(g(X)\mid S=1)$ affect the Disparate impact and the BER. Fairness should increase as the distance between the two distributions is small, which implies that $S$ does not affect the decision rule.  For this, in Figures \ref{fig:DI_W2} and \ref{fig:BER_W2}, we represent the evolution of the known criteria DI and BER with the Wasserstein distance $W_2^2(\mathcal{L}(g(X)\mid S=0),\mathcal{L}(g(X)\mid S=1))$, while the distributions $\mathcal{L}(X \mid S=0)$ and $\mathcal{L}(X \mid S=1)$ are being pushing forward onto their Wasserstein barycenter, according to the partial repair procedure called Geometric Repair \cite{FFMSV}, that moves each distribution part towards the Wasserstein barycenter in order to reduce the disparity between the groups.

\begin{figure}[h]
	\centering
	\begin{overpic}[scale=.5,unit=1mm]{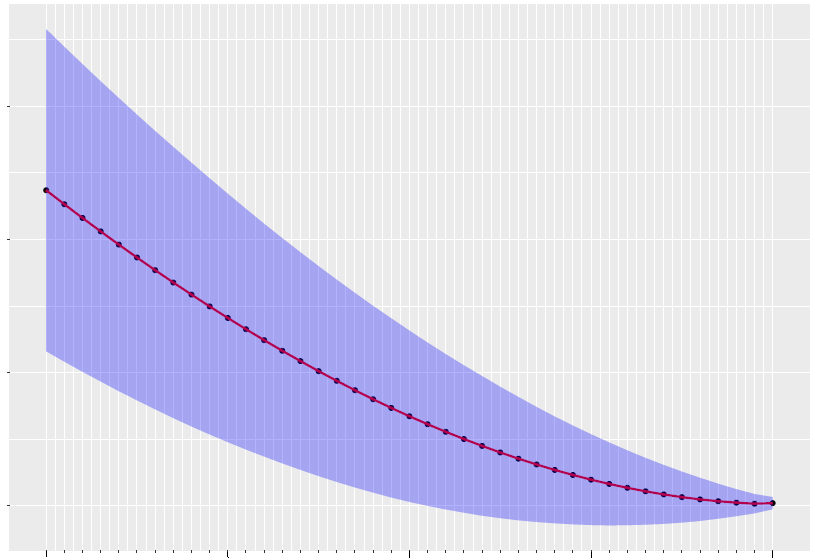}
		\put (-15,30){$(W_{n,2})^2$}
		
		\put (-3,5){$0$}
		\put (-8,22){$0.005$}
		\put (-8,38){$0.010$}
		\put (-8,54){$0.015$}
		
		\put (35,-10){Amount of repair}
		\put (5,-4){$0$}
		\put (25,-4){$0.25$}
		\put (48,-4){$0.5$}
		\put (70,-4){$0.75$}
		\put (94,-4){$1$}
		
	\end{overpic}
	\vspace{1cm}
	\caption{Confidence interval \eqref{asymptCI} for $W_{2}^2$}
	\label{fig:CI_W2}
\end{figure}

Figure \ref{fig:CI_W2} shows confidence intervals for the empirical quadratic transportation cost as the amount of repair increases. \vskip .1in

In Figure \ref{fig:DI_W2} we can see that the Disparate Impact decreases with the Wasserstein distance, and the desirable level equal to $0.8$ is attained when $W_2^2(\mathcal{L}(g(X)\mid S=0),\mathcal{L}(g(X)\mid S=1)) < 0.00225$. Note that 0.8 is a threshold chosen in many trials about unfair algorithmic treatments (see for instance in \cite{FFMSV} or \cite{ZVGG}). Moreover, Figure \ref{fig:BER_W2} confirms that the closer are both distributions in Wasserstein distance, the more unpredictable is the protected variable from the outcome of the regression. \vskip .1in
In conclusion although the distance in Total Variation and the Wasserstein distance are of very different nature, controlling the amount of fairness using the Wasserstein distance provides a control on the Disparate Impact. Moreover, it may be an alternative to the Balance Error and can provide a new control over the fairness of an algorithm. In this paper, we restricted ourselves to the logit classification but using a multidimensional version of the CLT as in \cite{del2017central}, we could provide a natural criterion of fairness directly on the observations $X$ by looking at the distance between $\mathcal{L}(X|S=1)$ and $\mathcal{L}(X|S=0)$.


\section*{Appendix}

In this Appendix we provide the proof of Theorem \ref{CLTLpCost}. Parts i) and ii) can be handled similarly. Hence, for the sake of simplicity we
focus on part i). The same techniques yield ii) with little extra effor. Throughout the section
we will assume that $U_1,\ldots, U_n$ are i.i.d. r.v.'s uniformly distributed on the interval $(0,1)$.
We write $A_n$ for the empirical distribution function on $U_1,\ldots,U_n$ and $\alpha_n(x)=\sqrt{n}(A_n(x)-x)$, $0\leq x\leq 1$
for the related empirical process. These $U_1,\ldots,U_n$ allow to represent any other 
i.i.d. sample $X_1,\ldots,X_n$ with d.f. $F$ by taking $X_i=F^{-1}(U_i)$. We use this construction in the sequel 
without further mention.

Given a distribution function $F$ we write $F_n$ for the empirical distribution function
based on the sample $F^{-1}(U_1),\ldots,F^{-1}(U_n)$ and $F_n^{-1}$ for the quantile inverse of $F_n$. Note that
$F_n^{-1}(t)=F^{-1}(A_n^{-1}(t))$. We fix a d.f. $G\in \mathcal{F}_{2p}(\mathbb{R})$ and define
\begin{equation}\label{Tnp}
T_{n,p}(F,G)=\sqrt{n}(\mathcal{W}_p^p(F_n,G)-E(\mathcal{W}_p^p(F_n,G))), \quad F\in \mathcal{F}_{2p}(\mathbb{R}).
\end{equation}
Similarly, using the notation in (\ref{GausspCost}) for $c_p$ and $\bar{c}_p$, we denote 
\begin{equation}\label{GaussianTransportProc}
T_p(F,G)=\int_0^1 \bar{c}_p(t;F,G) dW(t),
\end{equation}
where $\{W(t)\}_{0\leq t\leq 1}$ is a standard Brownian motion on $[0,1]$.
It follows from Lemma \ref{integrability} that $T_p(F,G)$ is a centered Gaussian r.v. with variance $\sigma^2(F,G)$
as in (\ref{VarF}).

\medskip 
We provide now some empirical counterparts of Lemma \ref{integrability}. First, a general variance bound for $T_{n,p}(F,G)$
and them, under more restrictive assumptions, an approximate continuity result for the trajectories of $T_{n,p}(\cdot,G)$.
The main ingredient in the proof is the Efron-Stein inequality for variances, namely, that if $Z=f(X_1,\ldots,X_n)$ with
$X_1,\ldots,X_n$ independent random variables, $(X_1',\ldots,X_n')$ is an independent copy of $(X_1,\ldots,X_n)$ and
$Z_i=f(X_1,\ldots,X_i',\ldots,X_n)$ then
$$\mbox{Var}(Z)\leq \sum_{i=1}^n E(Z-Z_i)_+^2.$$
We refer, for instance, to \cite{boucheron2013concentration} for further details.

\medskip

\begin{Proposition}\label{variancebound}
	If $F,G\in\mathcal{F}_{2p}$, $p>1$, then there exists a finite constant $C(F,G)$, depending only on $F$ and $G$ such that
	$$\mbox{\em Var}\big(T_{n,p}(F,G)\big)\leq {C(F,G)}, \quad n\geq 1.$$ 
	A valid choice of the constant is given by $C(F,G)=8p^2 \max(1,2^{2(p-1)})(C_1(F)+C_2(F,G))$ with
	$$C_1(F)=E\Big(|F^{-1}(U_1)-F^{-1}(U_2)|^2|F^{-1}(U_1)|^{2(p-1)}\Big)$$ 
	and
	$$C_2(F,G)=\Big(E\Big(|F^{-1}(U_1)-F^{-1}(U_2)|^{2p}\Big)\Big)^{1/p} \Big(E\Big(|G^{-1}(U_1)|^{2p}\Big)\Big)^{(p-1)/p}.
	$$
\end{Proposition}

\medskip
\noindent \textbf{Proof.} 
We recall that $F_n$ in equation (\ref{Tnp}) is the empirical 
distribution function based on the i.i.d. sample $X_i=F^{-1}(U_i)$, $i=1,\ldots,n$.
We set $Z=\mathcal{W}_p^p(F_n,G)$ and 
$Z'=\mathcal{W}_p^p(F_{n}',G)$, where $F_{n}'$ is the empirical distribution function based on 
the sample $X_1',X_2,\ldots,X_n$ and $X_1,X_1',X_2\ldots,X_n$ are i.i.d.. We write $X_{(1)}\leq \cdots \leq X_{(n)}$ 
for the ordered sample. Let us assume that $F$ is continuous. Now,
$Z=\sum_{i=1}^n  \int_{\frac {i-1} n}^{\frac i n} |X_{(i)}-G^{-1}(t)|^p dt  =\sum_{i=1}^n \int_{\frac{R_i-1} n}^{\frac{R_i}n} |X_i-G^{-1}(t)|^p dt$ with 
$R_i$ denoting the rank of $X_i$ within the sample $X_1,\ldots,X_n$. Continuity of $F$ ensures that a.s. there are no ties and
$(R_1,\ldots,R_n)$ is a random permutation of $\{1,\ldots,n\}$. Let us write $(R_1',\ldots,R_n')$ for the ranks in the sample
$X_1',X_2,\ldots,X_n$. Now, $Z$ is the minimal value of $E(|U-V|^p |X_1,\ldots,X_n, X_1')$ among random vectors $(U,V)$
which, conditionally given the $X_i$'s, have marginals $F_n$ and $G$. This shows that
$$ Z\leq \sum_{i=1}^n \int_{\frac{{R_i}'-1} n}^{\frac{{R_i}'} n} |X_i-G^{-1}(t)|^p dt$$
and, as a consequence, 
$$Z-Z'\leq \int_{\frac{{R_1}'-1} n}^{\frac{{R_1}'} n} \big[|X_1-G^{-1}(t)|^p-|X_1'-G^{-1}(t)|^p \big] dt.$$
Using the fact that
$||a+h|^{p}-|a|^p|\leq p\, h\,(|a+h|^{p-1}+|a|^{p-1})$ for $a\in\mathbb{R}, h>0, p>1$ and writing $d_p$ for the same
constants as in the proof of Lemma \ref{integrability}, we get that
\begin{eqnarray*}
	Z-Z'&\leq& p |X_1-X_1'|\int_{\frac{{R_1}'-1} n}^{\frac{{R_1}'} n} \big[|X_1-G^{-1}(t)|^{p-1}+|X_1'-G^{-1}(t)|^{p-1} \big] dt\\
	&\leq & p d_p |X_1-X_1'|\Big( 2\int_{\frac{{R_1}'-1} n}^{\frac{{R_1}'} n} |G^{-1}(t)|^{p-1}dt + \frac {|X_1|^{p-1}}{n}+\frac {|X_1'|^{p-1}}{n}\Big).
\end{eqnarray*}
Hence, 
$$E(Z-Z')_+^2\leq 8p^2 d_p^2 \Big[{\textstyle \frac 1 {n^2}}E \big(|X_1-X_1'|^2 |X_1|^{2p-2} \big)+
E \big(|X_1-X_1'|\int_{(R_1-1)/n}^{R_1/n} |G^{-1}(t)|^{p-1} dt  \big)^2\Big].$$
Under the assumption $F\in\mathcal{F}_{2p}$, $C_1(F):=E \big(|X_1-X_1'|^2 |X_1|^{2p-2} \big)$ is finite. To bound the last term we note that,
\begin{eqnarray*}
	E \Big(|X_1-X_1'|{\textstyle \int_{\frac{R_1-1}{n}}^{\frac {R_1} n} }|G^{-1}(t)|^{p-1} dt  \Big)^2&\leq &\big(E|X_1-X_1'|^{2p}\big)^{\frac 1 p} 
	\Big(E\Big({\textstyle \int_{\frac{R_1-1}{n}}^{\frac {R_1} n} }|G^{-1}(t)|^{p-1} dt\Big)^{\frac {2p}{p-1}}\Big)^{\frac{p-1}p}.
\end{eqnarray*}
Using again H\"older's inequality we see that
$$\Big({\textstyle \int_{\frac{j-1}{n}}^{\frac {j} n} } |G^{-1}(t)|^{p-1} dt\Big)^{\frac{2p}{p-1}}\leq {{n^{-\frac{p+1}{p-1}}}}
{ \int_{\frac{j-1}{n}}^{\frac {j} n} } |G^{-1}(t)|^{2p} dt$$
and, therefore,
\begin{eqnarray*}
	E\Big({\textstyle \int_{\frac{R_1-1}{n}}^{\frac {R_1} n} } |G^{-1}(t)|^{p-1} dt\Big)^{\frac {2p}{p-1}}&=&\frac 1 n\sum_{j= 1}^n 
	\Big({\textstyle \int_{\frac{j-1}{n}}^{\frac {j} n} }|G^{-1}(t)|^{p-1} dt\Big)^{\frac {2p}{p-1}}\\
	&\leq & {n^{-\frac{2p}{p-1}}} \sum_{j=1}^n {\textstyle \int_{\frac{j-1}{n}}^{\frac {j} n} } |G^{-1}(t)|^{2p} dt=\frac 1 
	{n^{\frac{2p}{p-1}}}\int_0^1 |G^{-1}(t)|^{2p} dt.
\end{eqnarray*}
As a consequence, 
$$E \Big(|X_1-X_1'|{\textstyle \int_{\frac{R_1-1}{n}}^{\frac {R_1} n} }|G^{-1}(t)|^{p-1} dt  \Big)^2\leq \frac {C_2(F,G)}{n^2},$$
with $C_2(F,G)=(E|X_1-X_1'|^{2p}\big)^{\frac 1 p} \big(\int_0^1 |G^{-1}(t)|^{2p} dt\big)^{\frac{p-1}{p}}<\infty$. Now the Efron-Stein inequality,
and the fact that $Z$ is a symmetric function of $X_1,\ldots,X_n$, which are i.i.d.
yields 
$$\mbox{Var}\big(\mathcal{W}_p^p(F_n,G)\big)\leq n  E(Z-Z')_+^2 \leq \frac{C(F,G)}n$$
with $C(F,G)=8p^2 d_p^2(C_1(F)+C_2(F,G))$. This yields the conclusion for continuous $F$. For general $F$ we take
continuous $F_m\in \mathcal{F}_{2p}(\mathbb{R})$ such that $\mathcal{W}_{2p}(F_m,F)\to 0$ as $m\to\infty$. A standard uniform
integrability argument shows that both $C(F_m,G)\to C(F,G)$ and $\mbox{Var}\big(T_{n,p}(F_m,G))\to
\mbox{Var}\big(T_{n,p}(F,G))$ as $m\to\infty$ and completes the proof.
\quad $\Box$

\medskip
An interesting consequence of Proposition \ref{variancebound} is that $T_{n,p}(F,G)$ can be approximated by
$T_{n,p}(F_M,$ $G_M)$ with $F_M, G_M$ being bounded support approximations of $F$ and $G$, respectively. We give details next.

\begin{Corollary}\label{truncation}
	Assume $F, G \in\mathcal{F}_{2p}(\mathbb{R})$ and $M>0$. Consider the distribution function $F_M$
	with quantile $F_M^{-1}(t)=\max(\min(F^{-1}(t),M),-M)$. Then there exist constants $C(M,F,G)$
	depending only on $M, F$ and $G$ such that
	$$\mbox{\em Var}(T_{n,p}(F,G)-T_{n,p}(F_M,G))\leq C(M,F,G),\quad n\geq 1$$
	and $C(M,F,G)\to 0$ as $M\to\infty$. Furthermore, if $G_M$ is 
	the distribution function with quantile $G_M^{-1}(t)=\max(\min(G^{-1}(t),M),-M)$
	then for every $\varepsilon>0$ there exist $M_0>0$ and $n_0$ such that
	$$\mbox{\em Var}(T_{n,p}(F,G)-T_{n,p}(F_M,G_M))\leq \varepsilon$$
	for each $M\geq M_0$ and $n\geq n_0$.
\end{Corollary}

\medskip
\noindent \textbf{Proof.} We write $\bar{F}_M$ for the distribution function with quantile
$\bar{F}_M^{-1}(t)=\min(F^{-1}(t),M)$. We will give a bound for 
$\mbox{Var}(T_{n,p}(F,G)-T_{n,p}(\bar{F}_M,G))$, with a similar argument for the left tail 
completing the proof. Now, observe that
\begin{eqnarray*}
	\lefteqn{T_{n,p}(F,G)-T_{n,p}(\bar{F}_M,G)=\sqrt{n}\Big[\int_0^1 |F^{-1}(A_n^{-1}(t))-G^{-1}(t)|^p dt-\int_0^1 |\bar{F}_M^{-1}(A_n^{-1}(t))-G^{-1}(t)|^p dt\Big]}
	\hspace*{2cm}\\
	&=&\sqrt{n}\Big[\int_{A_n^{-1}(t)>F(M)} |F^{-1}(A_n^{-1}(t))-G^{-1}(t)|^p dt-\int_{A_n^{-1}(t)>F(M)} |M-G^{-1}(t)|^p dt\Big].
\end{eqnarray*}
Note that the last expression does not depend on the values of $F^{-1}$ in the set $\{s\leq F(M)\}$. In particular,
if we write $\tilde{F}_M^{-1}(s)=F^{-1}(s)$, if $F^{-1}(s)>M$, $\tilde{F}_M^{-1}(s)=0$ otherwise, and  
$\hat{F}_M^{-1}(s)=M$, if $F^{-1}(s)>M$, $\hat{F}_M^{-1}(s)=0$ otherwise, then
$T_{n,p}(F,G)-T_{n,p}(\bar{F}_M,G)=T_{n,p}(\tilde{F}_M,G)-T_{n,p}(\hat{F}_M,G)$. As a consequence,
$$\mbox{Var}(T_{n,p}(F,G)-T_{n,p}(\bar{F}_M,G))\leq 2\mbox{Var}(T_{n,p}(\tilde{F}_M,G))+2\mbox{Var}(T_{n,p}(\hat{F}_M,G)).$$
It follows from Proposition \ref{variancebound} that
\begin{eqnarray*}
	\mbox{Var}(T_{n,p}(\tilde{F}_M,G))&\leq & \Big(\mu_{2p}(\tilde{F}_M)^{\frac{p-1}p} 
	+\mu_{2p}(G)^{\frac{p-1}p} \Big) \Big(E\Big(|\tilde{F}_M^{-1}(U_1)-\tilde{F}_M^{-1}(U_2)|^{2p}\Big)\Big)^{1/p}\\
	&\leq &
	\Big(\mu_{2p}({F})^{\frac{p-1}p} 
	+\mu_{2p}(G)^{\frac{p-1}p} \Big) \Big(E\Big(|\tilde{F}_M^{-1}(U_1)-\tilde{F}_M^{-1}(U_2)|^{2p}\Big)\Big)^{1/p},
\end{eqnarray*}
with $\mu_r(H)=\int_0^1 |H^{-1}(t)|^r dt.$ But $\tilde{F}_M^{-1}(U_1)-\tilde{F}_M^{-1}(U_2)$ vanishes
if $U_1\leq F(M)$ and $U_2\leq F(M)$. Hence,
\begin{eqnarray*}
	E\Big(|\tilde{F}_M^{-1}(U_1)-\tilde{F}_M^{-1}(U_2)|^{2p}\Big)\leq 2^{2p-1}
	\int_{((0,F(M))\times (0,F(M)) )^C} (|F^{-1}(s)|^{2p}+|F^{-1}(t)|^{2p})dsdt.
\end{eqnarray*}
By dominated convergence the last integral vanishes as $M\to\infty$. This proves the first claim. 
This allows to consider only the case of $F$ supported in $[-M,M]$ for the second claim. As before,
we show how to deal with the upper tail. Arguing as above, it suffices to bound the variance of $Z_M:=
\int_{G(M)}^1|F_n^{-1}(t)-G^{-1}(t)|^pdt$, with $F_n^{-1}(t)=F^{-1}(A_n^{-1}(t))$. We write
$X_i=F^{-1}(U_i)$, consider $X_1'$ and independent additional observation with law $F$ and argue as in the proof of 
Proposition \ref{variancebound}. We consider $Z_M'$, the version of $Z_M$ that we obtain replacing $X_1$ by
$X_1'$ in the sample and denote by $R_1$, $R_1'$ the ranks of $X_1$ and $X_1'$ in the samples. Now, if
$R_1\leq nG(M)$ and $R_1'\leq nG(M)$ then neither $X_1$ nor $X_1'$ enter in the expressions that define $Z_M$
and $Z_M'$, respectively, and, consequently, $Z_M-Z_M'=0$. Also, if $R_1'<R_1$ then $X_1'\leq X_1$ and (recall
that $X_1,\ldots,X_n,X_1'$ are upper bounded by $M$) replacing $X_1$ by $X_1'$ in the sample can only increase the transportation
cost, that is, $Z_M-Z_M'\leq 0$. Hence, if $Z_M-Z_M'>0$ then $R_1\leq R_1'$ and $R_1'>nG(M)$. If 
$R_1=R_1'$ then $Z_M-Z_M'\leq \int_{\frac{R_1'-1}{n}}^{\frac{R_1'}{n}} \Big| |X_1-G^{-1}(t)|^p-|X_1'-G^{-1}(t)|^p \Big| dt$.
If $R_1< R_1'$ then $X_1<X_1'$ and from the fact that $a<b<c<d$ implies $(d-b)^p+(c-a)^p\leq (d-a)^p+(b-c)^p$ we can see
that $Z_M-Z_M'\leq  \int_{\frac{R_1'-1}{n}}^{\frac{R_1'}{n}} \Big| |X_1-G^{-1}(t)|^p-|X_1'-G^{-1}(t)|^p \Big| dt$ as well.
Summarizing, we conclude that
$$Z_M-Z_M'\leq \int_{\frac{R_1'-1}{n}}^{\frac{R_1'}{n}} \Big| |X_1-G^{-1}(t)|^p-|X_1'-G^{-1}(t)|^p \Big| dt\,  I(R_1'>n G(M)).$$
We can now mimick the proof of Proposition \ref{variancebound} to see that
$$ E(Z_M-Z_M')_+^2\leq \frac{8p^2 d_p^2}{n^2} 
\Big(\mu_{2p}({F})^{\frac{p-1}p} 
+\mu_{2p}(G)^{\frac{p-1}p} \Big) \Big(E\Big(|X_1-X_1'|^{2p} I(R_1'>n G(M))\Big)\Big)^{1/p}.$$
Finally, we note that the probability that $R_1'$ exceeds $nG(M)$ is at most $1-G(M)+\frac 1 n$. This completes the proof. 
\quad $\Box$

\medskip
When $F$ and $G$ have bounded support and $G^{-1}$ is continuous it is possible to give
variance bounds for the increments of $T_{n,p}(\cdot,G)$. In view of Proposition
\ref{truncation} the assumption of bounded support does not mean a great loss in generality,
since slightly worse bounds can be obtained for the general case from this particular one.

\begin{Proposition}\label{CotaVarF1F2}
	If $F_1$, $F_2$ and $G$ are supported in $[-M,M]$ and $G^{-1}$ is continuous then
	there exists a sequence of constants $R_n(G,p,M)$, which depend on $G, p, M$ and $n$
	but not on $F_i$, $i=1,2$ such that $R_n(G,p,M)\to 0$ as $n\to\infty$ and
	$$\mbox{\em Var}(T_{n,p}(F_1,G)-T_{n,p}(F_2,G))\leq 3 \sigma_p^2(F_1,F_2;G)+M^2 R_n(G,p,M),$$
	with $\sigma_p^2(F_1,F_2;G)=E(T_p(F_1,G)-T_p(F_2,G))^2=\|\bar{c}_p(\cdot;F_1,G)- \bar{c}_p(\cdot;F_2,G) \|^2_{L_2(0,1)}$.
\end{Proposition}

\medskip
\noindent
\textbf{Proof.} We consider first a finitely supported $F$, 
concentrated on $x_1\leq\cdots\leq x_k$ with $F(x_j)=s_j$, $j=1,\ldots,k$. We have $s_k=1$
and set, for convenience, $s_0=0$. Then $\mathcal{W}_p^p(F,G)=\sum_{j=1}^k \int_{s_{j-1}}^{s_j} |x_{j}-G^{-1}(t)|^p dt$ and 
$\mathcal{W}_p^p(F_n,G)=\sum_{j=1}^k \int_{A_n(s_{j-1})}^{A_n(s_j)} |x_{j}-G^{-1}(t)|^p dt$. Hence,
$$\mathcal{W}_p^p(F,G)=\int_0^1 |x_k-G^{-1}(t)|^pdt-\sum_{j=1}^{k-1}\int_{0}^{s_j}\Big[ |x_{j+1}-G^{-1}(t)|^p -
|x_{j}-G^{-1}(t)|^p\Big] dt$$
and similarly for $\mathcal{W}_p^p(F_n,G)$, replacing $s_j$ with $A_n(s_j)$. Writing again $h_p(x)=|x|^p$
we have $|x_{j+1}-G^{-1}(t)|^p - |x_{j}-G^{-1}(t)|^p =\int_{x_j}^{x_{j+1}} h_p'(s-G^{-1}(t))ds$ and combining these 
last two facts, we obtain
$$T_n:=\sqrt{n}\Big(\mathcal{W}_p^p(F_n,G)-\mathcal{W}_p^p(F,G) \Big)=-\sqrt{n} \sum_{j=1}^{k-1}\int_{s_j}^{A_n(s_j)}\Big(
\int_{x_j}^{x_{j+1}} h_p'(s-G^{-1}(t))ds\Big) dt.$$
Next, we define 
$$\tilde{T}_n:=-\sqrt{n} \sum_{j=1}^{k-1}\int_{s_j}^{A_n(s_j)}\Big(
\int_{x_j}^{x_{j+1}} h_p'(s-G^{-1}(s_j))ds\Big) dt=- \sum_{j=1}^{k-1} \alpha_n(s_j)
\int_{x_j}^{x_{j+1}} h_p'(s-G^{-1}(s_j))ds$$
and observe that
\begin{equation}\label{cotaequicont}
|T_n-\tilde{T}_n|\leq \sqrt{n} \sum_{j=1}^{k-1}  \Big|\int_{s_j}^{A_n(s_j)}  \Big(
\int_{x_j}^{x_{j+1}} ( h_p'(s-G^{-1}(t))- h_p'(s-G^{-1}(s_j)))ds\Big) dt \Big|
\end{equation}
We consider now the continuity moduli
$$w_{G^{-1}}(\delta)=\sup_{x,y\in[0,1],|x-y|\leq \delta}|G^{-1}(x)-G^{-1}(y)|,$$
$$w_{p,M}(\varepsilon)=\sup_{x,y\in[-2M,2M],|x-y|\leq \varepsilon}|h_p'(x)(x)-h_p'(y)|.$$
The assumptions on $G^{-1}$ imply that it can be extended to a continuous function on $[0,1]$. Hence,
it is uniformly continuous and $w_{G^{-1}}(\delta)\to 0$ as $\delta\to 0$. Similarly,
$w_{p,M}(\varepsilon)\to 0$ as $\varepsilon\to 0$. Observe now that, for $t$ between $s_j$ and $A_n(s_j)$,
$|G^{-1}(t)-G^{-1}(s_j)|\leq w_{G^{-1}} (\|\alpha_n\|_\infty/\sqrt{n})$. Hence,
$$\int_{x_j}^{x_{j+1}} | h_p'(s-G^{-1}(t))- h_p'(s-G^{-1}(t))|ds\leq (x_{j+1}-x_j) w_{p,M}(  w_{G^{-1}} (\|\alpha_n\|_\infty/\sqrt{n})) $$
and, therefore, in view of (\ref{cotaequicont}),
\begin{eqnarray*}
	|T_n-\tilde{T}_n|&\leq &\sum_{j=1}^{k-1} (x_{j+1}-x_j) w_{p,M}(  w_{G^{-1}} (\|\alpha_n\|_\infty/\sqrt{n})) |\alpha_n(s_j)|\\
	&\leq & \|\alpha_n\|_\infty w_{p,M}(  w_{G^{-1}} (\|\alpha_n\|_\infty/\sqrt{n})) (x_k-x_1)\\
	&\leq & 2M \|\alpha_n\|_\infty w_{p,M}(  w_{G^{-1}} (\|\alpha_n\|_\infty/\sqrt{n})).
\end{eqnarray*}
Hence,
$$E(T_n-\tilde{T}_n)^2\leq M^2 \tilde{R}_n(G,p,M)$$
with $\tilde{R}_n(G,p,M)=4 E\Big[\|\alpha_n\|_\infty^2 w_{p,M}^2(  w_{G^{-1}} (\|\alpha_n\|_\infty/\sqrt{n})))\Big]$.
Uniform integrability of $\|\alpha_n\|_\infty^2$ and the fact that $w_{p,M}(  w_{G^{-1}} (\|\alpha_n\|_\infty/\sqrt{n}))$
is bounded and vanishes in probability ensure that $R_n(G,p,M)\to 0$ as $n\to\infty$.

Let us assume now that $F_1$ and $F_2$ are finitely supported as above and write $T_{n,i}$, $\tilde{T}_{n,i}$, 
$i=1,2$ for the corresponding versions of $T_n$ and $\tilde{T}_n$, respectively. Observe that there is no loss of generality in assumming
that $F_1$ and $F_2$ have a common support. Then
$$\mbox{Var}(T_{n,p}(F_1,G)-T_{n,p}(F_1,G))\leq E(T_{n,1}-T_{n,2})^2\leq 3 E(T_{n,1}-\tilde{T}_{n,1})^2+3
E(\tilde{T}_{n,1}-\tilde{T}_{n,2})^2+3E(T_{n,2}-\tilde{T}_{n,2})^2.$$
A simple covariance computation shows that $E(\tilde{T}_{n,1}-\tilde{T}_{n,2})^2=\sigma_p^2(F_1,F_2;G)$ and yields the conclusion.

For general $F_1$ and $F_2$ we take $F_{i,m}$, $i=1,2$, $m\geq 1$ with finite support (contained in $[-M,M]$) such that
$\mathcal{W}_{2p}(F_{i,m},F_i)\to 0$, $i=1,2$, and the bound follows by continuity. \quad $\Box$

\medskip

As a consequence of the variance bounds in Propositions \ref{variancebound} and \ref{CotaVarF1F2} and in Corollary
\ref{truncation} we can prove now the announced CLT for the empirical transportation cost.

\medskip
\noindent
\textbf{Proof of Theorem \ref{CLTLpCost}.} We will prove that
\begin{equation}\label{CLTLpCostWass}
\mathcal{W}_2\big(\mathcal{L}\big(T_{n,p}(F,G)\big), \mathcal{L}\big(T_{p}(F,G)\big)\big)\to 0
\end{equation}
As in the proof of Proposition \ref{truncation}, we assume first that $F$ is concentrated on 
$x_1\leq\cdots\leq x_k$ with $F(x_j)=s_j$, $j=1,\ldots,k$ and $F,G$ supported in $[-M,M]$. 
Then we have
$$T_n:=\sqrt{n}\Big(\mathcal{W}_p^p(F_n,G)-\mathcal{W}_p^p(F,G) \Big)=-\sqrt{n} \sum_{j=1}^{k-1}\int_{s_j}^{A_n(s_j)}\Big(
\int_{x_j}^{x_{j+1}} h_p'(s-G^{-1}(t))ds\Big) dt.$$
Continuity of $G^{-1}$ and the multivariate CLT imply that
$$\Big\{\sqrt{n}\int_{s_j}^{A_n(s_j)}\Big(
\int_{x_j}^{x_{j+1}} h_p'(s-G^{-1}(t))ds\Big) dt \Big\}_{j=1}^{k-1}\to_w
\Big\{ B(s_j)\int_{s_j}^{A_n(s_j)}\Big(
\int_{x_j}^{x_{j+1}} h_p'(s-G^{-1}(t))ds\Big) dt \Big\}_{j=1}^{k-1}$$
as $n\to\infty$, with $B(t)$ a Brownian bridge on $[0,1]$. 
Hence, using the trivial fact that $\sum_{j=0}^{k-1}B(s_j) c_j=-\sum_{j=0}^{k-1} d_j (B(s_{j+1})-B(s_j))$
if $d_0=c_0$ and $d_j=\sum_{l=0}^j c_l$, we conclude that
\begin{equation}\label{firstconv}
T_n\to_w T_p(F,G).
\end{equation}
We note that the assumptions on $F$ and $G$ guarantee that 
$$\Big|\sqrt{n}\int_{s_j}^{A_n(s_j)}\Big(
\int_{x_j}^{x_{j+1}} h_p'(s-G^{-1}(t))ds\Big) dt\Big| \leq K |\alpha_n(s_j)|$$
for some constant $K$. This shows that $T_n^2$ is uniformly integrable and, together with
(\ref{firstconv}) that $\mathcal{W}_2\Big(\mathcal{L}\Big(T_n\Big), \mathcal{L}\Big(T_p(F,G)\Big)\Big)\to 0$.
But this, in turn, yields convergence of moments of order 2 or smaller. In particular, we see that
$E(T_n)\to E(T_{p}(F,G))=0$, that is
\begin{equation}\label{centeringconstants}
\sqrt{n}\big(E(\mathcal{W}_p^p(F_n,G))-\mathcal{W}_p^p(F,G))\big)\to 0
\end{equation}
as $n\to\infty$. But (\ref{firstconv}) and (\ref{centeringconstants}) show that $T_{n,p}(F,G)\to_w T_p(F,G)$
and, again the uniform integrability, that 
$\mathcal{W}_2\Big(\mathcal{L}\Big(T_{n,p}(F,G)\Big), \mathcal{L}\Big(T_p(F,G)\Big)\Big)\to 0$.

\medskip
In a second step, we consider $F,G$ supported in $[-M,M]$,  with $G^{-1}$ continuous. We consider an approximating
sequence $F_m$ with finite support contained in $[-M,M]$ such that $\mathcal{W}_{2p}(F_m,F)\to 0$. Now, for
a fixed $\varepsilon>0$ we can, by Lemma \ref{integrability}, ensure that $\sigma_p^2(F_m,F,G)\leq \varepsilon^2$
for large $m$. For such $m$ we take $n_0$ large enough to guarantee that $R_n(G,p,M)\leq \varepsilon^2 /M^2$
and $\mathcal{W}_2\big(\mathcal{L}\big(T_{n,p}(F_m,G)\big), \mathcal{L}\big(T_p(F_m,$  $G)\big)\big)$ $\leq \varepsilon$
for $n\geq n_0$ (here $R_n(G,p,M)$ is as in Proposition \ref{CotaVarF1F2}). But then, for $n\geq n_0$,
\begin{eqnarray*}
	\lefteqn{\mathcal{W}_2\big(\mathcal{L}\big(T_{n,p}(F,G)\big), \mathcal{L}\big(T_{p}(F,G)\big)\big)\leq 
		\mathcal{W}_2\big(\mathcal{L}\big(T_{n,p}(F,G)\big), \mathcal{L}\big(T_{n,p}(F_m,G)\big)\big)}\hspace*{2.5cm}\\
	&+&\mathcal{W}_2\big(
	\mathcal{L}\big(T_{n,p}(F_m,G)\big), \mathcal{L}\big(T_p(F_m,G)\big)\big)
	+\mathcal{W}_2\big(\mathcal{L}\big(T_{p}(F_m,G)\big), \mathcal{L}\big(T_p(F,G)\big)\big)\\
	&\leq & 2\varepsilon + \varepsilon + \varepsilon=4\varepsilon,
\end{eqnarray*}
and we conclude that $\mathcal{W}_2\big(\mathcal{L}\big(T_{n,p}(F,G)\big), \mathcal{L}\big(T_{p}(F,G)\big)\big)\to 0$
as $n \to \infty$.

Finally, for $F,G \in\mathcal{F}_{2p}(\mathbb{R})$, $G^{-1}$ continuous we use Corollary \ref{truncation}. Note that $G_M^{-1}$
is also continuous. The already considered cases show that, for fixed $M$, 
$\mathcal{W}_2\big(\mathcal{L}\big(T_{n,p}(F_M,G_M)\big), \mathcal{L}\big(T_{p}(F_M,G_M)\big)\big)\to 0$ as $n\to\infty$.
Obviously, $\mathcal{W}_2\big(\mathcal{L}\big(T_{p}(F_M,G_M)\big), \mathcal{L}\big(T_{p}(F,G)\big)\big)\to 0$
as $M\to\infty$. Let us fix $\varepsilon>0$. We take $M_0$ and $n_0$ large enough to ensure that 
$\mathcal{W}_2\big(\mathcal{L}\big(T_{n,p}(F,G)\big),\mathcal{L}\big(T_{n,p}(F_M,G_M)\big)\big)\leq \epsilon$
if $M\geq M_0$ and $n\geq n_0$ and take $M\geq M_0$ large enough to guarantee 
$\mathcal{W}_2\big(\mathcal{L}\big(T_{p}(F_M,G_M)\big), \mathcal{L}\big(T_{p}(F,G)\big)\big)\leq \epsilon$. For this
choice of $M$ we take $n_1\geq n_0$ such that $\mathcal{W}_2\big(\mathcal{L}\big(T_{n,p}(F_M,G_M)\big), \mathcal{L}\big(T_{p}(F_M,G_M)\big)\big)
\leq \varepsilon$ for $n\geq n_1$. But then, arguing as above we see that 
$\mathcal{W}_2\big(\mathcal{L}\big(T_{n,p}(F,G)\big), \mathcal{L}\big(T_{p}(F,G)\big)\big)\leq 3 \varepsilon$ if $n\geq n_1$.
This completes the proof. \quad $\Box$

\medskip
\noindent \textbf{Proof of Proposition \ref{Proposition}.} As before, we give a proof for part \textit{(i)}. 
We will show first that under the given assumptions 
\begin{equation}\label{tightness}
\sqrt{n}(\mathcal{W}_p^p(F_n,G)-\mathcal{W}_p^p(F,G))\underset{w}\rightarrow N(0,\sigma^2)
\end{equation}
for some $\sigma^2\geq 0$.
For this goal we note that, by  assumption (\ref{tech1}),
$$\sqrt{n} \int_0^1 |F^{-1}-G^{-1}|^{p-2}|F_n^{-1}-F^{-1}|\leq \sqrt{n} (\mathcal{W}_p(F,G))^{{p-2}}
(\mathcal{W}_p(F_n,F))^{ 2}=o_P(1).$$
Similarly, we see that 
$$\sqrt{n} \int_0^1 |F_n^{-1}-G^{-1}|^{p-2}|F_n^{-1}-F^{-1}|=o_P(1).$$
A Taylor expansion and the fact that $|x|^{p-2}$ is a convex function imply that
\begin{eqnarray*}
	\lefteqn{\Big||F_n^{-1}-F^{-1}|^p-|F^{-1}-F^{-1}|^p-(F_n^{-1}-F^{-1})h'_p(F^{-1}-G^{-1}) \Big|}\hspace*{2cm}\\
	&\leq & C 
	(F_n^{-1}-F^{-1})^2 \Big(|F^{-1}-G^{-1}|^{p-2}+|F_n^{-1}-G^{-1}|^{p-2}  \Big).
\end{eqnarray*}
This bound and the above estimates yield that
$$\sqrt{n}\big(\mathcal{W}_p^p(F_n,G)- \mathcal{W}_p^p(F,G)-\int_0^1 (F_n^{-1}-F^{-1})h_p'(F^{-1}-G^{-1})\big)
=o_P(1).$$
Hence, we focus on the analysis of $\sqrt{n}\int_0^1 (F_n^{-1}-F^{-1})h_p'(F^{-1}-G^{-1})$. 
The moment assumptions on
$F$ and $G$ (see, e.g., Lemma 3.3 in \cite{alvarez2011uniqueness}) allow to replace  $\int_0^1 $ by $\int_{\frac 1 n}^{1-\frac 1 n}$ without
modyfing the asymptotic behavior of the resulting r.v.. Also, by Lemma 2.3 in \cite{del2005asymptotics} and
assumptions (\ref{smoothnessplus}) and (\ref{tech2}), we can replace $\sqrt{n}(F_n^{-1}-F^{-1})$ in the integral by the weighted
uniform
quantile process, $u_n/f(F^{-1}(\cdot))$, where $u_n(t)=\sqrt{n}(A_n^{-1}(t)-t)$. Therefore, to prove (\ref{tightness}) it suffices
to prove convergence of
$$\textstyle \int_{1/n}^{1-1/n} \frac{u_n}{f(F^{-1}(\cdot))} h_p'(F^{-1}-G^{-1}).$$ 
But now, Theorem 4.2 in \cite{del2005asymptotics}, assumptions (\ref{tech2}) and (\ref{tech3}) and the fact that $h_p'(F^{-1}-G^{-1})$
yield the result.

Now, from (\ref{tightness}) and Theorem \ref{CLTLpCost}  
we conclude that $\sqrt{n}(E\mathcal{W}_p^p(F_n,G)-\mathcal{W}_p^p(F,G))$
must be bounded. This in turn 
yields moment convergence (up to order two; recall the proof of Theorem \ref{CLTLpCost})
of $\sqrt{n} (\mathcal{W}_p^p(F_n,G)-\mathcal{W}_p^p(F,G))$. But since the limiting distribution of 
$\sqrt{n} (\mathcal{W}_p^p(F_n,G)-\mathcal{W}_p^p(F,G))$ is, as noted above, centered, we must
have 
$$\sqrt{n}(E\mathcal{W}_p^p(F_n,G)-\mathcal{W}_p^p(F,G))\to 0.$$
This concludes the proof. \hfill $\Box$


\begin{thebibliography}{99}

\bibitem{ajtai1984optimal}
Mikl{\'o}s Ajtai, J{\'a}nos Koml{\'o}s, and G{\'a}bor Tusn{\'a}dy.
\newblock On optimal matchings.
\newblock {\em Combinatorica}, 4(4):259--264, 1984.

\bibitem{alvarez2011uniqueness}
PC~{\'A}lvarez-Esteban, E~Del~Barrio, JA~Cuesta-Albertos, C~Matr{\'a}n, et~al.
\newblock Uniqueness and approximate computation of optimal incomplete
  transportation plans.
\newblock In {\em Annales de l'Institut Henri Poincar{\'e}, Probabilit{\'e}s et
  Statistiques}, volume~47, pages 358--375. Institut Henri Poincar{\'e}, 2011.

\bibitem{berk2017fairness}
Richard Berk, Hoda Heidari, Shahin Jabbari, Michael Kearns, and Aaron Roth.
\newblock Fairness in criminal justice risk assessments: the state of the art.
\newblock {\em arXiv preprint arXiv:1703.09207}, 2017.

\bibitem{myfairtest}
Philippe Besse, Eustasio Del~Barrio, Paula Gordaliza, and Jean-Michel Loubes.
\newblock Confidence intervals for testing disparate impact in fair learning.
\newblock {\em ArXiv e-prints 1807.06362}, 2018.

\bibitem{bobkov2014one}
Sergey Bobkov and Michel Ledoux.
\newblock One-dimensional empirical measures, order statistics and kantorovich
  transport distances.
\newblock {\em preprint}, 2014.

\bibitem{boucheron2013concentration}
St{\'e}phane Boucheron, G{\'a}bor Lugosi, and Pascal Massart.
\newblock {\em Concentration inequalities: A nonasymptotic theory of
  independence}.
\newblock Oxford university press, 2013.

\bibitem{2017arXiv170300056C}
A.~{Chouldechova}.
\newblock {Fair prediction with disparate impact: A study of bias in recidivism
  prediction instruments}.
\newblock {\em ArXiv e-prints}, February 2017.

\bibitem{del2018obtaining}
Eustasio del Barrio, Fabrice Gamboa, Paula Gordaliza, and Jean-Michel Loubes.
\newblock Obtaining fairness using optimal transport theory.
\newblock {\em arXiv preprint arXiv:1806.03195}, 2018.

\bibitem{del1999central}
Eustasio del Barrio, Evarist Gin{\'e}, and Carlos Matr{\'a}n.
\newblock Central limit theorems for the wasserstein distance between the
  empirical and the true distributions.
\newblock {\em Annals of Probability}, pages 1009--1071, 1999.

\bibitem{del2005asymptotics}
Eustasio Del~Barrio, Evarist Gin{\'e}, Frederic Utzet, et~al.
\newblock Asymptotics for l2 functionals of the empirical quantile process,
  with applications to tests of fit based on weighted wasserstein distances.
\newblock {\em Bernoulli}, 11(1):131--189, 2005.

\bibitem{del2017central}
Eustasio Del~Barrio and Jean-Michel Loubes.
\newblock Central limit theorem for empirical transportation cost in general
  dimension.
\newblock {\em arXiv preprint arXiv:1705.01299}, 2017.

\bibitem{dobric1995asymptotics}
V~Dobri{\'c} and Joseph~E Yukich.
\newblock Asymptotics for transportation cost in high dimensions.
\newblock {\em Journal of Theoretical Probability}, 8(1):97--118, 1995.

\bibitem{FFMSV}
Michael Feldman, Sorelle~A Friedler, John Moeller, Carlos Scheidegger, and
  Suresh Venkatasubramanian.
\newblock Certifying and removing disparate impact.
\newblock In {\em Proceedings of the 21th ACM SIGKDD International Conference
  on Knowledge Discovery and Data Mining}, pages 259--268. ACM, 2015.

\bibitem{fournier2015rate}
Nicolas Fournier and Arnaud Guillin.
\newblock On the rate of convergence in wasserstein distance of the empirical
  measure.
\newblock {\em Probability Theory and Related Fields}, 162(3-4):707--738, 2015.

\bibitem{2018arXiv180204422F}
S.~A. {Friedler}, C.~{Scheidegger}, S.~{Venkatasubramanian}, S.~{Choudhary},
  E.~P. {Hamilton}, and D.~{Roth}.
\newblock {A comparative study of fairness-enhancing interventions in machine
  learning}.
\newblock {\em ArXiv e-prints}, February 2018.

\bibitem{DBLP:journals/corr/abs-1712-07924}
Philipp Hacker and Emil Wiedemann.
\newblock A continuous framework for fairness.
\newblock {\em CoRR}, abs/1712.07924, 2017.

\bibitem{pedreschi2012study}
Dino Pedreschi, Salvatore Ruggieri, and Franco Turini.
\newblock A study of top-k measures for discrimination discovery.
\newblock In {\em Proceedings of the 27th Annual ACM Symposium on Applied
  Computing}, pages 126--131. ACM, 2012.

\bibitem{ramdas2017wasserstein}
Aaditya Ramdas, Nicol{\'a}s~Garc{\'\i}a Trillos, and Marco Cuturi.
\newblock On wasserstein two-sample testing and related families of
  nonparametric tests.
\newblock {\em Entropy}, 19(2):47, 2017.

\bibitem{romei_ruggieri_2014}
Andrea Romei and Salvatore Ruggieri.
\newblock A multidisciplinary survey on discrimination analysis.
\newblock {\em The Knowledge Engineering Review}, 29(5):582?638, 2014.

\bibitem{sommerfeld2018inference}
Max Sommerfeld and Axel Munk.
\newblock Inference for empirical wasserstein distances on finite spaces.
\newblock {\em Journal of the Royal Statistical Society: Series B (Statistical
  Methodology)}, 80(1):219--238, 2018.

\bibitem{Talagrand}
M~Talagrand and JE~Yukich.
\newblock The integrability of the square exponential transportation cost.
\newblock {\em The Annals of Applied Probability}, pages 1100--1111, 1993.

\bibitem{tameling2017empirical}
Carla Tameling, Max Sommerfeld, and Axel Munk.
\newblock Empirical optimal transport on countable metric spaces:
  Distributional limits and statistical applications.
\newblock {\em arXiv preprint arXiv:1707.00973}, 2017.

\bibitem{Villani2003}
C{\'e}dric Villani.
\newblock {\em Topics in optimal transportation}.
\newblock Number~58. American Mathematical Soc., 2003.

\bibitem{ZVGG}
Muhammad~Bilal Zafar, Isabel Valera, Manuel Gomez~Rodriguez, and Krishna~P
  Gummadi.
\newblock Fairness beyond disparate treatment \& disparate impact: Learning
  classification without disparate mistreatment.
\newblock In {\em Proceedings of the 26th International Conference on World
  Wide Web}, pages 1171--1180. International World Wide Web Conferences
  Steering Committee, 2017.

\end{thebibliography}
\end{document}